\documentclass[a4paper,11pt]{article}
\usepackage{amsmath,amsfonts,amssymb,amsthm}

\usepackage{graphicx}
\usepackage[all]{xy}
\usepackage{xcolor}

\def\spn{\bigskip\par\noindent}
\def\mpn{\medskip\par\noindent}
\def\pn{\par\noindent}

\def\hkr{\hookrightarrow}

\def\CE{\mathcal{E}}
\def\CF{\mathcal{F}}
\def\CO{\mathcal{O}}
\def\CP{\mathcal{P}}

\def\CS{\mathcal{S}}

\def\Hom{\mathrm{Hom}}

\def\zero{\{0\}}

\def\un{\mathbf{1}}
\def\Ind{\mathrm{Ind}}
\def\Inf{\mathrm{Inf}}

\def\Res{\mathrm{Res}}
\def\Iso{\mathrm{Iso}}
\def\Out{\mathrm{Out}}

\def\Irr{\mathrm{Irr}}

\def\End{\mathrm{End}}

\def\sfc{\mathsf{c}}
\def\SR{\mathsf{R}}

\definecolor{vert}{rgb}{0.12, 0.6, 0.17}

\def\sur{\overline}

\def\Proj{{\mathsf{Proj}}}
\def\smp{\smallskip\par}
\def\dom{\backslash}

\def\calC{\mathcal{C}}

\def\calF{\mathcal{F}}
\def\calP{\mathcal{P}}
\def\calL{\mathcal{L}}

\def\ZZ{\mathbb{Z}}
\def\NN{\mathbb{N}}
\newcommand{\pf}{{\flushleft\bf Proof: }}
\def\endpf{~\leaders\hbox to 1em{\hss\  \hss}\hfill~\raisebox{.5ex}{\framebox[1ex]{}}\smp}

\newcounter{nonce}[section]

\newenvironment{enonce}[1]{\pagebreak[2]\refstepcounter{nonce}\vskip 2ex plus 1ex minus 1ex\par\noindent{{\bf #1~\thesection.\arabic{nonce}:}}\begin{it}}{\end{it}\vskip 2ex plus 2ex minus 2ex\pn}
\newenvironment{enonce*}[1]{\pagebreak[2]\vskip 2ex plus 1ex minus 1ex\par\noindent{{\bf #1:}}\begin{it}}{\end{it}\vskip 2ex plus 2ex minus 2ex\pn}
\newenvironment{rem}[1]{\pagebreak[2]\refstepcounter{nonce}\vskip 2ex plus 1ex minus 1ex\par\noindent{{\bf #1~\thesection.\arabic{nonce}:}}}{\vskip 2ex plus 2ex minus 2ex\pn}
\newenvironment{moneq}{\refstepcounter{nonce}\begin{equation}}{\end{equation}}
\newcommand{\genk}[1]{\langle k_{#1}\rangle}
\newcommand{\ora}[1]{\overrightarrow{#1}}
\newcommand{\ola}[1]{\overleftarrow{#1}}
\begin{document}
\centerline{\Large\bf The diagonal $p$-permutation functor $kR_k$}\vspace{2ex}\par
\centerline{\bf Serge Bouc}\vspace{2ex}
\begin{abstract} Let $k$ be an algebraically closed field of positive characteristic $p$. We describe the full lattice of subfunctors of the diagonal $p$-permutation functor $kR_k$ obtained by $k$-linear extension from the functor $R_k$ of linear representations over $k$. This leads to the description of the ``composition factors'' $S_P$ of $kR_k$, which are parametrized by finite $p$-groups (up to isomorphism), and of the evaluations of these particular simple diagonal $p$-permutation functors over $k$.\end{abstract}

\section{Introduction} Let $k$ be an algebraically closed field of positive characteristic $p$, and $\SR$ be a commutative and unital ring. We consider\footnote{mainly in the case $\SR=k$, but a few results deal with the case of an arbitrary $\SR$.} the assignment $G\mapsto \SR R_k(G)$ sending a finite group $G$ to $\SR R_k(G)=\SR\otimes_\ZZ R_k(G)$, where $R_k(G)$ is the Grothendieck group of finite dimensional $kG$-modules. When $G$ and $H$ are finite groups, if $M$ is a $(kH,kG)$-bimodule which is projective as a right $kG$-module, tensoring with $M$ over $kG$ preserves exact sequences, hence induces a well defined group homomorphism $R_k(M):R_k(G)\to R_k(H)$, and an $\SR$-linear map $\SR R_k(M):\SR R_k(G)\to \SR R_k(H)$, both simply  denoted $V\mapsto M\otimes_{kG}V$. In particular, one checks easily that this endows $\SR R_k$ with a structure of {\em diagonal $p$-permutation functor} over $\SR$ (see \cite{bouc-yilmaz}, \cite{functorialequivalence}, \cite{diagonal-functors-in-char-p}). Our main goal is Theorem~\ref{main}, describing the full lattice of subfunctors of $kR_k$ (ordered by inclusion of subfunctors).\par
As a biproduct of this description, we get a ``composition series'' of $kR_k$, with simple subquotients $S_P$ indexed by finite $p$-groups $P$ up to isomorphism. From this, the evaluations of these particular simple functors can be computed, giving a new proof of the last corollary of~\cite{diagonal-functors-in-char-p}.\par
To be more precise, we first need some notation and some definitions:
\begin{enonce}{Notation} Let $\mathcal{P}$ be the set of isomorphism classes of finite $p$-groups. For a finite $p$-group $P$, we denote by $[P]$ its isomorphism class. For a subset $\calC$ of $\calP$, we say (abusively) that $P$ belongs to $\calC$, and we write $P\in\calC$, if $[P]\in\calC$.\par
For a finite group $G$, we denote by $k_G$ the trivial $kG$-module $k$, and by $\genk{G}$ the subfunctor of $kR_k$ generated by $k_G\in kR_k(G)$. We observe that $\genk{G}$ only depends on the isomorphism class of $G$. \par
For a finite group $G$, we denote by $\CE_\SR(G)$ the {\em essential algebra} of $G$ over~$k$, defined by
$$\CE_\SR(G)=\SR T^\Delta(G,G)/\sum_{|H|<|G|}\SR T^\Delta(G,H)\circ \SR T^\Delta(H,G),$$
where $T^\Delta(H,G)$ denotes the Grothendieck group of diagonal $p$-permutation $(kH,kG)$-bimodules (see \cite{bouc-yilmaz}, \cite{functorialequivalence}, \cite{diagonal-functors-in-char-p} for details).
\end{enonce}
\begin{enonce}{Definition} We say that a subset $\calC$ of $\calP$ is {\em closed} if for any finite $p$-groups $P$ and $Q$
$$Q\leq P\;\hbox{and}\;[P]\in\calC\;\implies\;[Q]\in\calC.$$ 
\end{enonce}
{\flushleft In} other words, if a $p$-group is isomorphic to a subgroup of a $p$-group in $\calC$, then it belongs to $\calC$.\par
Our main theorem is the following:
\begin{enonce}{Theorem} \label{main}Let $\calL$ denote the lattice of subfunctors of $kR_k$, ordered by inclusion of subfunctors. Let $\calF$ denote the lattice of closed subsets of $\calP$, ordered by inclusion of subsets. Then the maps
\begin{align*}
\Psi: F\in\calL\mapsto \Psi(F)&=\big\{[P]\in\calP\mid k_P\in F(P)\big\}\in\calF,\,\hbox{and}\\
\Theta: \mathcal{C}\in\calF\mapsto \Theta(\calC)&=\sum_{[P]\in\calC}\genk{P}\in\calL,
\end{align*}
\nopagebreak are well defined isomorphisms of lattices, inverse to each other.
\end{enonce}
\section{Some evaluations of diagonal $p$-permutation functors.} We start with a result of independent interest on diagonal $p$-permutation functors. It gives a way to compute the evaluation of a diagonal $p$-permutation functor at a direct product of a $p'$-group and a $p$-group. This requires the following notation:
\begin{enonce}{Notation} Let $L$ be a $p'$-group and $Q$ be a $p$-group, and set $H=L\times Q$. For a (finite dimensional) $kL$-module $V$, let $\ora{V}=V\otimes_k kQ$, wiewed as a $(kH,kQ)$-bimodule for the action
$$\forall (l,x)\in H=L\times Q, \,\forall y,z\in Q, \,\forall v\in V,\;(l,x)\cdot (v\otimes z)\cdot y=lv\otimes xzy.$$
Similarly, let $V^*$ denote the $k$-dual of $V$, viewed as a right $kL$-module, and $\ola{V}=V^*\otimes_k kQ$, viewed as a $(kQ,kH)$-bimodule for the action
$$\forall (l,x)\in H=L\times Q, \,\forall y,z\in Q, \,\forall \alpha\in V^*,\;y\cdot(\alpha\otimes z)\cdot(l,x)=\alpha l\otimes yzx.$$
\end{enonce}
{\flushleft Let} also $\Irr_k(L)$ denote a set of representatives of isomorphism classes of irreducible $kH$-modules. With this notation:
\begin{enonce}{Proposition}\label{LxQ}
\begin{enumerate}
\item The bimodules $\ora{V}$ and $\ola{V}$ are diagonal $p$-permutation bimodules.
\item If $V$ and $W$ are irreducible $kL$-modules, then $\ola{W}\otimes_{kH}\ora{V}$ is isomorphic to the identity $(kQ,kQ)$-bimodule $kQ$, if $W$ and $V$ are isomorphic, and it is zero otherwise.
\item The direct sum $\mathop{\oplus}_{V}\limits \big(\ora{V}\otimes_{kQ}\ola{V}\big)$, for $V\in\Irr_k(L)$, is isomorphic to the identity $(kH,kH)$-bimodule $kH$.
\end{enumerate}
\end{enonce}
\pf 1. Indeed, the group $H\times Q$ has a normal Sylow subgroup $S=Q\times Q$, and the restrictions to $S$ of $\ora{V}$ and $\ola{V}$ are both isomorphic to a direct sum of $\dim_kV$ copies of the identity $(kQ,kQ)$-bimodule $kQ$.\mpn
2. One checks easily that
$$\ola{W}\otimes_{kH}\ora{V}=(W^*\otimes_k kQ)\otimes_{k(L\times Q)}(V\otimes_k kQ)\cong (W^*\otimes_{kL}V)\otimes_k(kQ\otimes_{kQ}kQ)$$
as $(kQ,kQ)$-bimodules. Assertion 2 follows, since $W^*\otimes_{kL}V=0$ unless $V$ and $W$ are isomorphic, and since $V^*\otimes_{kL}V\cong k$ as $k$ is algebraically closed. Moreover $kQ\otimes_{kQ}kQ\cong kQ$.\mpn
3. Similarly
\begin{align*}
\ora{V}\otimes_{kQ}\ola{V}&=(V\otimes_k kQ)\otimes_{kQ}(V^*\otimes_k kQ)\cong(V\otimes_k V^*)\otimes_k(kQ\otimes_{kQ}kQ)\\
&\cong \End_k(V)\otimes_kkQ
\end{align*}
as $(kH,kH)$-bimodules. Now the direct sum $\mathop{\oplus}_{V\in\Irr_k(L)}\limits\End_k(V)$ is the decomposition of the semisimple group algebra $kL$ as a direct product of its Wedderburn components. It follows that $\mathop{\oplus}_{V}\limits \big(\ora{V}\otimes_{kQ}\ola{V}\big)\cong kL\otimes kQ\cong k(L\times Q)$ is isomorphic to the identity $(kH,kH)$-bimodule.\endpf
In the following corollary, we view $\ora{V}$ and $\ola{V}$ as elements of $T^\Delta(H,Q)$ and $T^\Delta(Q,H)$, respectively.
\begin{enonce}{Corollary} \label{F(LxQ)}\label{essential}Let $\SR$ be a commutative (unital) ring, and $F$ be a diagonal $p$-permutation functor over $\SR$. Let $L$ be a $p'$-group, and $Q$ be a $p$-group. Then:
\begin{enumerate}
\item $F(L\times Q)\cong \mathop{\bigoplus}_{V\in \Irr_k(L)}\limits F(Q)$. More precisely, the maps
$$\varphi\in F(L\times Q)\mapsto \mathop{\bigoplus}_{V\in \Irr_k(L)}F(\ola{V})(\varphi)\in \mathop{\bigoplus}_{V\in \Irr_k(L)}F(Q)$$
and 
$$\psi=\mathop{\bigoplus}_{V\in \Irr_k(L)}\psi_V\in \mathop{\bigoplus}_{V\in \Irr_k(L)}F(Q)\mapsto \sum_{V\in \Irr_k(L)}F(\ora{V})(\psi_V)\in F(L\times Q)$$
are isomorphisms of $\SR$-modules, inverse to each other.
\item If $L$ is non-trivial, then the essential algebra $\CE_\SR(L\times Q)$ is zero.
\end{enumerate}
\end{enonce}
\pf 1. Indeed, setting $H=L\times Q$ as above, for irreducible $kL$-modules $V$ and $W$, we have that
$$F(\ola{W})\circ F(\ora{V})=F(\ola{W}\otimes_{kH}\ora{V}),\;\hbox{and}\;F(\ora{V})\circ F(\ola{V})=F(\ora{V}\otimes_{kQ}\ola{V}),$$
so the result follows from Assertions 2 and 3 of Proposition~\ref{LxQ}.\mpn
2. Assertion 3 of Proposition~\ref{LxQ} tells us that the identity bimodule $kH$ is equal to a sum of morphisms which factor through $Q$, and $|Q|<|H|$ if $L$ is non-trivial. Assertion 2 of Corollary~\ref{F(LxQ)} follows.
\endpf
\begin{rem}{Remark} \begin{enumerate}
\item In short, Corollary~\ref{F(LxQ)} says that
$$F(L\times Q)\cong R_k(L)\otimes_\ZZ F(Q).$$
\item The (non-)vanishing of $\CE_\SR(H)$ for an arbitrary finite group $H$ is studied in detail in Section~3 of~\cite{diagonal-functors-in-char-p}, under additional conditions on $\SR$. The proof given here is much simpler and explicit for the case of a direct product $H=L\times Q$ of a $p'$-group $L$ and a $p$-group $Q$, and doesn't require any additional assumption on $\SR$.
\end{enumerate} 
\end{rem}
We now consider a version of Brauer's induction theorem relative to the prime~$p$:
\begin{enonce}{Theorem} Let $G$ be a finite group. Then the cokernel of the induction map
$$\mathop{\oplus}_H\limits\Ind_H^G:\bigoplus_H R_k(H)\to R_k(G)$$
where $H$ runs through Brauer $p$-elementary subgroups of $G$, i.e. subgroups of the form $H=P\times C$, where $P$ is a $p$-group and $C$ is a cyclic $p'$-group, is finite and of order prime to $p$.
\end{enonce}
\pf  Let $(K,\mathcal{O},k)$ be a $p$-modular system, such that $K$ is big enough for~$G$. We have a commutative diagram
$$\xymatrix{
\mathop{\bigoplus}_H\limits R_K(H)\ar[r]^-{\mathop{\oplus}_H\limits\Ind_H^G}\ar[d]_-{\mathop{\oplus}_H\limits d_H}&R_K(G)\ar[d]^-{d_G}\\
\mathop{\bigoplus}_H\limits R_k(H)\ar[r]_-{\mathop{\oplus}_H\limits\Ind_H^G}&R_k(G)\\
}
$$
where the vertical arrows $d_H$ and $d_G$ are the respective decomposition maps. Let $C_K$ (resp. $C_k$) denote the cokernel of the top (resp. bottom) horizontal map. By Theorem~12.28 of~\cite{groupesfinis}, the groups $C_K$ is finite, and of order prime to~$p$. By Theorem 16.33 of~\cite{groupesfinis}, the vertical arrows are surjective. So $d_G$ induces a surjective group homomorphism $C_K\twoheadrightarrow C_k$, hence $C_k$ is finite and of order prime to $p$.\endpf
\begin{enonce}{Corollary} \label{Brauer elementary}Let $G$ be a finite group. Then
$$kR_k(G)=\sum_{H}\Ind_{H}^G\,kR_k(H),$$
where $H$ runs through the Brauer $p$-elementary subgroups of $G$.
\end{enonce}
\begin{enonce}{Corollary} \label{ideal} Let $F$ be a subfunctor of $kR_k$, and $G$ be a finite group. Then $F(G)$ is an ideal of the algebra $kR_k(G)$. 
\end{enonce}
\pf By Corollary~\ref{Brauer elementary}, there is a set $\CS$ of Brauer elementary subgroups of $G$ and elements $w_H\in kR_k(H)$, for $H\in\CS$, such that
$$k_G=\sum_{H\in\CS}\Ind_H^Gw_H.$$
now let $u\in F(G)$. Then in $kR_k(G)$, we have
$$u=k_G\cdot u=\sum_{H\in\CS}\Ind_H^G(w_H\cdot\Res_H^Gu).$$
Since $F$ is a subfunctor of $kR_k$, we know that $\Res_H^GF(G)\subseteq F(H)$ and $\Ind_H^GF(H)\subseteq F(G)$. So it is enough to prove that $F(H)$ is an ideal of $kR_k(H)$, for each $H$ in $\CS$.\par
Now each $H\in \CS$ is in particular of the form $H=L\times Q$, where $L$ is a (cyclic) $p'$-group, and $Q$ is a $p$-group. In particular $R_k(H)=\Inf_L^{L\times Q}R_k(L)$. Let $V$ and $W$ be $kL$-modules, and let $\widetilde{V}$ denote the $(kL,kL)$-bimodule $\Ind_{\Delta(L)}^{L\times L}V$. Then $\widetilde{V}\otimes_k kQ$ is a $(kH,kH)$-diagonal $p$-permutation bimodule. Moreover, one can check easily that $\widetilde{V}\otimes_{kL}W\cong V\otimes_k W$ as $kL$-modules, and it follows that there is an isomorphism of $kH$-modules
$$(\widetilde{V}\otimes_{k}kQ)\otimes_{k(L\times Q)}\Inf_{L}^{L\times Q}W\cong \Inf_{L}^{L\times Q}(V\otimes_kW).$$
Written differently, this reads
$$kR_k\big(\widetilde{V}\otimes_{k}kQ\big)\big(\Inf_{L}^{L\times Q}W\big)=\Inf_{L}^{L\times Q}(V\otimes_kW).$$
In the algebra $kR_k(H)=\Inf_L^{L\times Q}kR_k(L)$, the right hand side is nothing but the product $V\cdot W$, and the left hand side is given by the action of the diagonal $p$-permutation bimodule $\widetilde{V}\otimes_k kQ$ on $\Inf_L^{L\times Q}W$. It follows more generally that if $\varphi\in F(H)$, then the product $V\cdot \varphi$ is obtained from $\varphi$ by applying $\widetilde{V}\otimes_k kQ\in T^\Delta(H,H)$, so $V\cdot\varphi\in F(H)$. Thus $kR_k(H)\cdot F(H)\subseteq F(H)$, as was to be shown.
\endpf
\begin{enonce}{Notation} For finite groups $G$ and $H$, we write $H\hkr G$ if $H$ is isomorphic to a subgroup of $G$.
\end{enonce}
\begin{enonce}{Lemma} \label{inclusions}Let $P$ and $Q$ be finite $p$-groups. Then 
$$\genk{P}(Q)=\left\{\begin{array}{ll}k &\hbox{if }\;Q\hookrightarrow P\\0&\hbox{otherwise}.\end{array}\right.$$
In particular $\genk{Q}\leq\genk{P}$ if and only if $Q\hookrightarrow P$.
\end{enonce}
\pf By definition of $\genk{P}$, we have $\genk{P}(Q)=kT^\Delta(Q,P)(k_P)$. Moreover since $P$ and $Q$ are $p$-groups, the group $T^\Delta(Q,P)$ is equal to the Burnside group $B^\Delta(Q,P)$ of left-right free $(Q,P)$-bisets, so $kT^\Delta(Q,P)=kB^\Delta(Q,P)$. Now if $X$ is a diagonal subgroup of $Q\times P$, the set of (right-)orbits of $P$ on the $(Q,P)$-biset $(Q\times P)/X$ is a $Q$-set isomorphic to $Q/p_1(X)$, where $p_1(X)$ is the first projection of $X$, so we have
$$
k\big((Q\times P)/X\big)\otimes_{kP}k\cong k\big(Q/p_1(X)\big).
$$
As an element of $R_k(Q)$, this is equal to $|Q:p_1(X)|k_Q$, so it is equal to 0 in $kR_k(Q)$ unless $p_1(X)=Q$, and then it is equal to $k_Q$. But saying that there exists a diagonal subgroup $X\leq Q\times P$ such that $p_1(X)=Q$ amounts to saying that $Q$ is isomorphic to a subgroup of $P$. \par
Since $\genk{P}$ is a subfunctor of $kR_k$, and since $kR_k(Q)$ is one dimensional, generated by $k_Q$, it follows that $\genk{P}(Q)$ is non zero if and only if $Q$ is isomorphic to a subgroup of $P$, and this occurs if and only if $k_Q\in\genk{P}(Q)$, that is $\genk{Q}\leq\genk{P}$. This completes the proof.\endpf
\section{The subfunctors $\genk{P}$}
The following result is a major step in the proof of Theorem~\ref{main}:
\begin{enonce}{Theorem} \label{subfunctor at G}Let $F$ be a subfunctor of $kR_k$. Then for any finite group $G$
$$F(G)=\sum_{\substack{P\in\calP,\,P\leq G\\k_P\in F(P)}}\genk{P}(G).$$
\end{enonce}
\pf Let $G$ be a finite group. Saying that $k_P\in F(P)$ amounts to saying that $\genk{P}\leq F$, so $\mathop{\displaystyle\sum}_{\substack{P\in\calP, P\leq G\\k_P\in F(P)}}\limits\genk{P}(G)\leq F(G)$. We have to prove the reverse inclusion.\par
By Theorem~\ref{Brauer elementary}, there is a set $\CS$ of Brauer $p$-elementary subgroups of~$G$, and elements $w_{H}\in kR_k(H)$, for $H\in \CS$, such that
$$k_G=\sum_{H\in\CS}\Ind_{H}^Gw_{H}$$
in $kR_k(G)$.
Let $u\in F(G)$. Then
\begin{moneq}\label{Frobenius}
u=k_G\otimes_ku=\sum_{H\in\CS}\Ind_{H}^G\big(w_H\otimes_k\Res_{H}^Gu\big).
\end{moneq}
{\flushleft Now} $\Res_{H}^Gu\in F(H)$, for each $H\in\CS$, since $F$ is a subfunctor of $kR_k$. Then $w_H\otimes_k\Res_{H}^Gu\in F(H)$ also, by Lemma~\ref{ideal}, since each $H\in\CS$ is a product $L_H\times Q_H$, where $L_H$ is a (cyclic) $p'$-group and $Q_H$ is a $p$-group. It follows that $u\in \sum_{H\in\CS}\limits\Ind_H^GF(H)$ for any $u\in F(G)$, so
$$F(G)=\sum_{H\in\CS}\Ind_H^GF(H).$$
Now by Proposition~\ref{LxQ}, for each $H\in\CS$, the identity $(kH,kH)$-bimodule $kH$ splits as
$$kH\cong \mathop{\oplus}_{V\in\Irr_k(L_H)}\big(\ora{V}\otimes_{kQ_H}\ola{V}\big).$$
It follows that
$$F(H)=\sum_{V\in\Irr_k(L_H)}F(\ora{V})\big(F(Q_H)\big).$$
Now $F(Q_H)\subseteq kR_k(Q_H)=k$, since $Q_H$ is a $p$-group. So $F(Q_H)$ is either zero or $k$. It is non zero exactly when $k_{Q_H}\in F(Q_H)$. So
$$F(H)=\sum_{V\in\Irr_k(L_H)}F(\ora{V})\big(F(Q_H)\big)=\genk{Q_H}(H),$$
and this is non-zero if and only if $k_{Q_H}\in F(Q_H)$. It follows that
$$F(G)=\sum_{\substack{H\in\CS\\k_{Q_H}\in F(Q_H)}}\Ind_H^G\genk{Q_H}(H)\subseteq \sum_{\substack{P\in\calP,\,P\leq G\\k_P\in F(P)}}\genk{P}(G),$$
as was to be shown.\endpf

\begin{enonce}{Corollary} \label{subfunctor}Let $F$ be a subfunctor of $kR_k$. Then
$$F=\sum_{\substack{P\in \CP\\k_P\in F(P)}}\limits\genk{P}.$$
\end{enonce}
\begin{enonce}{Corollary} \label{Sylow}Let $G$ be a finite group, and $P$ be a Sylow $p$-subgroup of $G$. Then:
\begin{enumerate}
\item $kR_k(G)=\genk{P}(G)$.
\item $\genk{G}=\genk{P}$.
\end{enumerate}
\end{enonce}
\pf 1. Indeed $kR_k(G)=\mathop{\displaystyle\sum}_{Q\in\calP,\,Q\leq G}\limits\genk{Q}(G)$, by Theorem~\ref{subfunctor at G} applied to the subfunctor $F=kR_k$. But if $Q\in\calP$ and $Q\leq G$, then $Q\hkr P$, hence $\genk{Q}\leq\genk{P}$ by Lemma~\ref{inclusions}. Hence $\genk{Q}(G)\leq\genk{P}(G)$, and $kR_k(G)=\genk{P}(G)$. \spn
2. Since $k_G\in\genk{G}(G)$, we have $k_G\in\genk{P}(G)$ by Assertion 1, that is $\genk{G}\leq\genk{P}$. Conversely $k_P=\Res_P^Gk_G$, so $k_P\in \genk{G}(P)$, i.e. $\genk{P}\leq\genk{G}$. Hence $\genk{G}=\genk{P}$.
\endpf
\section{Proof of Theorem~\ref{main}}
We first recall the statement:
\begin{enonce*}{Theorem} Let $\calL$ denote the lattice of subfunctors of $kR_k$, ordered by inclusion of subfunctors. Let $\calF$ denote the lattice of closed subsets of $\calP$, ordered by inclusion of subsets. Then the maps
\begin{align*}
\Psi: F\in\calL\mapsto \Psi(F)&=\big\{[P]\in\calP\mid k_P\in F(P)\big\}\in\calF,\,\hbox{and}\\
\Theta: \mathcal{C}\in\calF\mapsto \Theta(\calC)&=\sum_{[P]\in\calC}\genk{P}\in\calL,
\end{align*}
\nopagebreak are well defined isomorphisms of lattices, inverse to each other.
\end{enonce*} 
\pf We first check that $\Psi$ is well defined, i.e. that $\Psi(F)$ is a closed subset of $\calP$, for any $F\in\calL$. This follows from Lemma~\ref{inclusions}: Saying that $P\in\Psi(F)$ amounts to saying that $\genk{P}\leq F$. Now if $Q\hookrightarrow P\in\Psi(F)$, Lemma~\ref{inclusions} shows that $\genk{Q}\leq\genk{P}\leq F$, so $Q\in\Psi(F)$, as was to be shown.\par
The maps $\Psi$ and $\Theta$ are clearly maps of posets. Moreover Corollary~\ref{subfunctor} shows that $\Theta\circ\Psi(F)=F$, for any $F\in\calL$. Conversely, let $\calC$ be a closed subset of $\calP$. Then
$$\Psi\circ\Theta(\calC)=\big\{[P]\in\calP\mid k_P\in \sum_{Q\in\calC}\genk{Q}(P)\big\}.$$
So clearly $\calC\subseteq \Psi\circ\Theta(\calC)$, since $k_Q\in \genk{Q}(Q)$. Conversely, if $[P]\in \Psi\circ\Theta(\calC)$, that is if $k_P$ belongs to $\sum_{Q\in\calC}\limits\genk{Q}(P)$, there is some $Q\in\calC$ such that $\genk{Q}(P)$ is non zero. Then $P\hookrightarrow Q$, by Lemma~\ref{inclusions}, so $P\in\calC$ since $\calC$ is closed. So $\Psi\circ\Theta(\calC)=\calC$, which completes the proof of Theorem~\ref{main}.\endpf
\begin{enonce}{Corollary} \begin{enumerate}
\item The lattice $\calL$ is completely distributive.
\item Let $F$ be a subfunctor of $kR_k$. The following are equivalent in the lattice~$\cal L$:
\begin{enumerate}
\item $F$ is completely join irreducible.
\item $F$ is completely prime.
\item There is a finite $p$-group $P$ such that $F=\genk{P}$.
\end{enumerate}
Moreover in (c), the $p$-group $P$ is unique up to isomorphism.
\end{enumerate} 
\end{enonce}
\pf 1. Indeed, the lattice $\calL$ is isomorphic to $\CF$, which is clearly completely distributive: its join is union of closed subsets, and its meet is intersection. Now if $\calC$ and $(\calC_i)_{i\in I}$ are closed subsets of $\CP$, then $\calC\cap \mathop{\cup}_{i\in I}\limits\calC_i=\mathop{\cup}_{i\in I}\limits(\calC\cap \calC_i)$.\spn
2. Clearly (b) implies (a). Moreover $(c)$ implies $(b)$: if $P$ is a $p$-group, then by Lemma~\ref{inclusions}
$$\Psi\big(\genk{P}\big)=\{Q\in\CP\mid Q\hkr P\}.$$
Hence if $(\calC_i)_{i\in I}$ is a set of closed subsets of $\CP$, then $\Psi\big(\genk{P}\big)\subseteq \mathop{\cup}_{i\in I}\limits \calC_i$ if and only if there exists some $i\in I$ such that $P\in\calC_i$, i.e. $\Psi\big(\genk{P}\big)\subseteq \calC_i$. So $\Psi\big(\genk{P}\big)$ is completely prime in $\calF$, and $\genk{P}$ is completely prime in $\calL$.\par
Finally (a) implies (c): if $F$ is completely join irreducible in $\calL$, since $F=\sum_{k_P\in F(P)}\limits\genk{P}$, there is a $p$-group $P$ such that $F=\genk{P}$. Moreover $P$ is unique up to isomorphism, by Lemma~\ref{inclusions}.
\endpf
\begin{rem}{Remark} There are join irreducible elements of $\calL$, or equivalently of $\calF$, which are {\em not} completely join irreducible: Let for example $\calC$ be the subset of $\CP$ consisting of cyclic $p$-groups. Then $\calC$ is closed, and one checks easily that $\calC$ is join irreducible in $\calF$, but not completely join irreducible. 
\end{rem}
\section{Some simple diagonal $p$-permutation functors}
\begin{enonce}{Proposition} \label{JP}Let $P$ be a finite $p$-group. Then $\genk{P}$ has a unique (proper) maximal subfunctor $J_P$, defined for a finite group $G$ by
$$J_P(G)=\big\{u\in \genk{P}(G)\mid kT^\Delta(P,G)(u)=0\big\}.$$
Moreover $J_P=\mathop{{\displaystyle\sum}}_{\substack{Q\hookrightarrow P\\Q\ncong P}}\limits\genk{Q}$.
\end{enonce}
\pf First it is easy to check that the assignment $G\mapsto J_P(G)$ defines a subfunctor of $kR_k$, hence of $\genk{P}$. Moreover $J_P(P)=0$. Now Lemma~\ref{inclusions} implies that $\genk{P}(P)=k$, so $J_P$ is a proper subfunctor of $\genk{P}$. If $F$ is a subfunctor of $\genk{P}$, there are two possibilities: Either $F(P)=k$, and then $k_P\in F(P)$, so $F=\genk{P}$. Or $F(P)=0$, and then for any finite group $G$, we have $kT^\Delta(P,G)\big(F(G)\big)\leq F(P)=0$, that is $F\leq J_P$. \par
For the last assertion, denote by $].,P]$ the subset of $\calP$ consisting of $p$-groups isomorphic to a subgroup of $P$, and by $].,P[$ the subset of $].,P]$ consisting of $p$-groups isomorphic to a {\em proper} subgroup of $P$. Then clearly $].,P]\in\mathcal F$, and $\Theta(].,P])=\genk{P}$, by Lemma~\ref{inclusions}. Also $].,P[\in\calF$, and $\Theta\big(].,P[\big)=\mathop{{\displaystyle\sum}}_{\substack{Q\hookrightarrow P\\Q\ncong P}}\limits\genk{Q}$, by definition of $\Theta$. Now$].,P[$ is clearly the unique maximal proper closed subset of $].,P]$, so $\Theta\big(].,P[\big)=J_P$, by Theorem~\ref{main}. This completes the proof.\endpf
\begin{enonce}{Notation} For a finite $p$-group $P$, we denote by 
$$S_P=\genk{P}/J_P$$
the unique simple quotient of $\genk{P}$.
\end{enonce}
\begin{enonce}{Lemma} \label{SP}Let $G$ be a finite group, and $P$ be a finite $p$-group.\begin{enumerate}
\item If $S_P(G)\neq 0$, then $P\hookrightarrow G$.
\item If $Q$ is a finite $p$-group, then 
$$S_P(Q)=\left\{\begin{array}{ll}k&\hbox{if $Q\cong P$}\\0&\hbox{otherwise.}\end{array}\right.$$
In particular $S_P\cong S_Q$ if and only if $P\cong Q$.
\end{enumerate}
\end{enonce}
\pf 1. By Theorem~\ref{subfunctor at G}, we have 
$$\genk{P}(G)=\sum_{\substack{Q\in\calP,\,Q\leq G\\k_Q\in \genk{P}(Q)}}\genk{Q}(G).$$
Now $k_Q\in \genk{P}(Q)$ if and only if $\genk{Q}\leq\genk{P}$, i.e. $Q\hkr P$, by Lemma~\ref{inclusions}. So if $P$ is not isomorphic to a subgroup of $G$, then $Q\leq G$ and $Q\hkr P$ implies that $Q$ is isomorphic to a proper subgroup of $P$, and then $\genk{Q}(G)\leq J_P(G)$, by Proposition~\ref{JP}. It follows that $\genk{P}(G)\leq J_P(G)$, so $S_P(G)=0$.\spn
2. Indeed, if $S_P(Q)\neq 0$, then in particular $\genk{P}(Q)\neq 0$, so $Q\hookrightarrow P$ by Lemma~\ref{inclusions}. But also $P\hkr Q$ by Assertion 1. So $S_P(Q)=0$ unless $Q\cong P$. Moreover $S_P(P)=\genk{P}(P)/J_P(P)=k/\{0\}\cong k$, which completes the proof.\endpf
\begin{enonce}{Theorem} \label{composition}\begin{enumerate}
\item Let $F_2<F_1$ be subfunctors of $kR_k$ such that $F_1/F_2$ is a simple functor. Then there exists a (unique, up to isomorphism) finite $p$-group $P$ such that $F_1/F_2\cong S_P$.
\item There exists a filtration
$$0=F_0< F_1<\ldots <F_n< F_{n+1}<\ldots$$
of $kR_k$ by subfunctors $F_i$, for $i\in\NN$, such that:
\begin{enumerate}
\item $\mathop{\bigcup}_{i=0}^\infty\limits F_i=kR_k$.
\item For $i>0$, the functor $F_i/F_{i-1}$ is simple, isomorphic to $S_{P_i}$, for a finite $p$-group $P_i$.
\item For every finite $p$-group $P$, there is exactly one integer $i>0$ such that $P_i\cong P$.
\end{enumerate}
\end{enumerate}
\end{enonce}
\pf 1. Set $\calC_i=\Psi(F_i)$, for $i\in\{1,2\}$. Then $\calC_1$ and $\calC_2$ are closed subsets of of $\calP$, and $\calC_2\subset\calC_1$. Since $F_1/F_2$ is simple, any subfunctor $F$ of $kR_k$ such that $F_2\leq F\leq F_1$ is equal either to $F_1$ or $F_2$. Then any closed subset $\calC$ of $\calP$ such that $\calC_2\subseteq \calC\subseteq \calC_1$ is equal either to $\calC_1$ or $\calC_2$. If $P\in\calC_1-\calC_2$, then $\calC_2\cup ].,P]$ is closed, different from $\calC_2$, and contained in $\calC_1$. So $\calC_2\cup ].,P]=\calC_1$. Now if $P'\in\calC_1-\calC_2$, then $P'\in].,P]$, and $P\in ].,P']$, by exchanging the roles of $P$ and $P'$. It follows that $P\cong P'$, so $\calC_1-\calC_2=\{[P]\}$. Now
$$F_1=\Theta(\calC_1)=\Theta(\calC_2)+\genk{P}=F_2+\genk{P}.$$
It follows that $F_1/F_2\cong \genk{P}/(\genk{P}\cap F_2)$ is a simple quotient of $\genk{P}$, so $F_1/F_2\cong S_P$. The uniqueness of $P$ (up to isomorphism) with this property follows from Lemma~\ref{SP}.\spn
2. Choose an enumeration $P_1, P_2,\ldots, P_n,\ldots$ of $\cal P$ such that for any indices $i$ and $j$, $P_i\hookrightarrow P_j$ implies $i\leq j$. This can be achieved starting with {$P_1=\nolinebreak\un$}, $P_2=C_p$, and then enumerating all the $p$-groups of order $p^2$, then the groups of order $p^3$, and so on. With such a numbering, set $\calC_0=\emptyset$ and $\calC_n=\cup_{i\leq n}].,P_i]$ for $n>0$, and then set $F_n=\Theta(\calC_n)$ for $n\in\NN$. In particular $F_0=0$.\par
Since $\calC_n=\calC_{n-1}\cup ].,P_n]$ for $n>0$, the sequence $(\calC_n)_{n\in \NN}$ is increasing. Moreover $P_n\in \calC_n-\calC_{n-1}$, for otherwise $P_n\in ].,P_i]$ for some $i<n$, meaning that $P_n\hookrightarrow P_i$, which implies $n\leq i$, a contradiction. So the sequence $(\calC_n)_{n\in \NN}$ is stricly increasing, and its union is the whole of $\calP$. In other words, the sequence $(F_n)_{n\in \NN}$ is stricly increasing and its union is equal to $kR_k$, which proves~(a).\par
Since $\calC_i=\calC_{i-1}\cup ].,P_i]$ for $i>0$, and since $].,P_i[\subseteq \calC_{i-1}$ by our choice of the numbering of $\calP$, it follows that $\calC_i=\calC_{i-1}\cup\{P_i\}$, and then $F_i/F_{i-1}\cong S_{P_i}$ as in the proof of Assertion 1. This proves (b). Now (c) is clear, since any finite $p$-group $P$ appears exactly once in our enumeration.\endpf
\begin{rem}{Remark} Theorem~\ref{composition} shows that $kR_k$ admits a ``composition series'', where the ``composition factors'' are the simple functors $S_P$, for $P\in\calP$, and each simple functor $S_P$ appears exactly once. The quote signs in the previous sentence indicate that one should be careful with the notions of composition factors and composition series for diagonal $p$-permutation functors.\vspace{-2ex}
\end{rem} 
\begin{rem}{Remark} \label{other parametrization}Proposition~\ref{JP} and Lemma~\ref{SP} show that $P$ is a minimal group for~$S_P$: More precisely $S_P(P)$ is one dimensional, and $S_P(G)=\zero$ for any group~$G$ of order smaller than $|P|$. Moreover $P$ is unique (up to isomorphism) with this property. In addition $S_P(P)$ is generated by the image of the trivial module $k_P$, and in particular the group $\Out(P)$ of outer automorphisms of $P$ acts trivially on $S_P(P)$. These two facts show that with the notation of Theorem 5.25 of~\cite{diagonal-functors-in-char-p}, the functor $S_P$ is isomorphic to the simple functor $\mathsf{S}_{P,1,k}$.
\end{rem}
\section{The simple functor $S_1$}\label{S un}
We consider first the case where the $p$-group $P$ is trivial\footnote{The content of the present section~\ref{S un} is essentially the same as Subsection 6.5 of~\cite{diagonal-functors-in-char-p}. We include it here for the reader's convenience.}. In this case $S_P=S_1=\genk{1}\leq kR_k$, so for a finite group $G$
$$S_1(G)=k T^\Delta(G,1)(k_1).$$
But $k T^\Delta(G,1)$ is the group $kP_k(G)$ of projective $kG$-modules. So $S_1(G)=\genk{1}(G)$ is equal to the image of the map $k\sfc^G:kP_k(G)\to kR_k(G)$ linearly extending the Cartan map $\sfc^G:P_kG)\to R_k(G)$ sending a projective $kG$-module to the sum of its composition factors. \par
\begin{rem}{Remark} \label{socle}The functor $S_1$ is the only simple subfunctor of $kR_k$: indeed, any non-empty closed subset of $\CP$ must contain the trivial group, so any non-zero subfunctor of $kR_k$ must contain $\genk{1}=S_1$.\par
As there are non-empty closed subsets of $\CP$ different from $\{1\}$, this  shows that $S_1$ is a proper subfunctor of $kR_k$, so $kR_k$ is not a semisimple diagonal $p$-permutation functor. It follows that in contrast to~\cite{functorialequivalence}, the category of diagonal $p$-permutation functors over a field $k$ {\em of characteristic $p$} is {\em not semisimple}.
\end{rem}
{\flushleft We} choose a $p$-modular system $(K,\CO,k)$, and we assume that $K$ is big enough for the group $G$. If $S$ is a simple $kG$-module, we denote by $\mathsf{P}_S$ its projective cover, and by $\Phi_S:G_{p'}\to \CO$ the modular character of $\mathsf{P}_S$, where $G_{p'}$ is the set of $p$-regular elements of $G$. If $v=\sum_{S\in \Irr_k(G)}\limits \omega_S\mathsf{P}_S$, where $\omega_S\in \CO$, is an element of $\CO\Proj(kG)$, we denote by $\Phi_v:\CO\Proj(kG)\to \CO$ the map $\sum_{S\in\Irr_k(G)}\limits \omega_S\Phi_S:G_{p'}\to \CO$, and we  call $\Phi_v$ the modular character of $v$. \par
Then for a simple $kG$-module $T$, the multiplicity of $S$ as a composition factor of $\mathsf{P}_T$ is equal to the Cartan coefficient
$$\sfc^G_{T,S}=\dim_k\Hom_{kG}(\mathsf{P}_T,\mathsf{P}_S)=\frac{1}{|G|}\sum_{x\in G_{p'}}\Phi_T(x)\Phi_S(x^{-1}).$$
In order to describe the image of the map $k\sfc^G$, we want to evaluate the image of this integer under the projection map $\rho:\CO\to k$. For this, we denote by $[G_{p'}]$ a set of representatives of conjugacy classes of $G_{p'}$, and we observe that in the field~$K$, we have
\begin{align}\refstepcounter{nonce}
\sfc^G_{S,T}&=\frac{1}{|G|}\sum_{x\in [G_{p'}]}\frac{|G|}{|C_G(x)|}\Phi_T(x)\Phi_S(x^{-1})\nonumber\\
&=\sum_{x\in [G_{p'}]}\frac{1}{|C_G(x)|}\frac{\Phi_T(x)}{|C_G(x)|_{p}}\frac{\Phi_S(x^{-1})}{|C_G(x)|_{p}}|C_G(x)|_{p}^2\nonumber\\
&=\sum_{x\in [G_{p'}]}\frac{\big(\Phi_T(x)/|C_G(x)|_{p}\big)\big(\Phi_S(x^{-1})/|C_G(x)|_{p}\big)}{|C_G(x)_{p'}|}|C_G(x)|_{p}.\label{in O}
\end{align}
But since $\Phi_S$ and $\Phi_T$ are characters of projective $kG$-modules (and since $C_G(x)=C_G(x^{-1}))$, the quotients $\Phi_T(x)/|C_G(x)|_{p}$ and $\Phi_S(x^{-1})/|C_G(x)|_{p}$ are in $\CO$, so 
$$\forall x\in[G_{p'}],\;\frac{\big(\Phi_T(x)/|C_G(x)|_{p}\big)\big(\Phi_S(x^{-1})/|C_G(x)|_{p}\big)}{|C_G(x)_{p'}|}\in\CO.$$
Now it follows from~\ref{in O} that
\begin{moneq}\label{rho}\rho(\sfc^G_{T,S})=\sum_{x\in [G_0]}\rho\left(\frac{\Phi_T(x)\Phi_S(x^{-1})}{|C_G(x)|}\right),
\end{moneq}
{\flushleft where} $[G_0]$ is a set of representatives of conjugacy classes of the set $G_0$ of elements {\em defect zero} of $G$, i.e. the set of $p$-regular elements $x$ such that $C_G(x)$ is a $p'$-group.
\begin{enonce}{Notation} \label{notation Gamma}For $x\in G_0$, we set
$$\Gamma_{G,\,x}=\sum_{S\in\Irr(kG)}\frac{\Phi_S(x^{-1})}{|C_G(x)|}S\in\CO R_k(G),$$
where $\Irr(kG)$ is a set of representatives of isomorphism classes of simple $kG$-modules. We also set
$$\gamma_{G,\,x}=\sum_{S\in\Irr(kG)}\rho\left(\frac{\Phi_S(x^{-1})}{|C_G(x)|}\right)S\in kR_k(G),$$
\end{enonce}
\begin{rem}{Remark}\label{class}
We note that $\Gamma_{G,\,x}$ and $\gamma_{G,\,x}$ only depend on the conjugacy class of $x$ in $G$, that is $\Gamma_{G,\,x}=\Gamma_{G,\,x^g}$ and $\gamma_{G,\,x}=\gamma_{G,\,x^g}$ for $g\in G$.
\end{rem}
{\flushleft By} Theorem 6.3.2 of~\cite{brauer-nesbitt} (see also Theorem 6.3.2 of~\cite{nagao-tsushima}), the elementary divisors of the Cartan matrix of $G$ are equal to $|C_G(x)|_p$, for $x\in [G_{p'}]$. It follows that the rank of the Cartan matrix modulo $p$, is equal to the number of conjugacy classes of elements of defect~0 of $G$, i.e. the cardinality of $[G_0]$. The following can be viewed as an explicit form of this result:

\begin{enonce}{Proposition}\label{basis}\begin{enumerate}
\item Let $T$ be a simple $kG$-module. Then, in $kR_k(G)$,
$$k\sfc^G(\mathsf{P}_T)=\sum_{x\in[G_0]}\rho\big(\Phi_T(x)\big)\,\gamma_{G,\,x}.$$
\item The elements $\gamma_{G,\,x}$, for $x\in[G_0]$, form a basis of $S_1(G)\leq kR_k(G)$.
\end{enumerate}
\end{enonce}
\pf Throughout the proof, we simply write $\gamma_x$ instead of $\gamma_{G,\,x}$.\mpn
1. By~\ref{rho}, we have
\begin{align*}
k\sfc^G(\mathsf{P}_T)&=\sum_{S\in\Irr(kG)}\rho(\sfc^G_{T,S})S=\sum_{S\in\Irr(kG)}\sum_{x\in [G_0]}\rho\left(\frac{\Phi_T(x)\Phi_S(x^{-1})}{|C_G(x)|}\right)S\\
&=\sum_{x\in [G_0]}\sum_{S\in\Irr(kG)}\rho\left(\frac{\Phi_T(x)\Phi_S(x^{-1})}{|C_G(x)|}\right)S\\
&=\sum_{x\in [G_0]}\rho\big(\Phi_T(x)\big)\sum_{S\in\Irr(kG)}\rho\left(\frac{\Phi_S(x^{-1})}{|C_G(x)|}\right)S\\
&=\sum_{x\in[G_0]}\rho\big(\Phi_T(x)\big)\,\gamma_x.
\end{align*}
2. We first prove that $\gamma_x$ lies in the image of $k\sfc^G$, for any $x\in G_0$. So let $x\in G_0$, and $1_x:\langle x\rangle \to \CO$ be the map with value 1 at $x$ and~0 elsewhere. Then $|x|1_x=\sum_{\zeta}\zeta(x^{-1})\zeta$, where $\zeta$ runs through the simple $k\langle x\rangle$-modules, i.e. the group homomorphisms $\langle x\rangle\to \CO^\times$, is an element of $\CO P_k(\langle x\rangle)=\CO R_k(\langle x\rangle)$. Let $v_x=\Ind_{\langle x\rangle}^G(|x|1_x)$. Then $v_x\in \CO P_k(G)$, and its modular character evaluated at $g\in G$ is equal to
\begin{align}\refstepcounter{nonce}
\Phi_{v_x}(g)&=\frac{1}{|x|}\sum_{\substack{h\in G\\g^h\in\langle x\rangle}}\Phi_{|x|1_x}(g^h)\nonumber\\
\label{character}&=\frac{1}{|x|}\sum_{\substack{h\in G\\g^h=x}}|x|=\left\{\begin{array}{ll}|C_G(x)|&\hbox{if}\;g=_Gx\\0&\hbox{otherwise},\end{array}\right.
\end{align}
where $g=_Gx$ means that $g$ is conjugate to $x$ in $G$. Now from Assertion 1, we get that 
\begin{moneq}\label{v viz gamma}
k\sfc^G(v_x)=\sum_{y\in[G_0]}\rho\big(\Phi_{v_x}(y)\big)\gamma_y=|C_G(x)|\gamma_x,
\end{moneq}
so $\gamma_x$ is in the image of $k\sfc^G$, since $|C_G(x)|\neq 0$ in $k$.\par
Now by Assertion 1, the elements $\gamma_x$, for $x\in[G_0]$, generate the image of $k\sfc^G$, i.e. $S_1(G)$. They are moreover linearly independent: Suppose indeed that some linear combination $\sum_{x\in [G_0]}\limits\lambda_x\gamma_x$, where $\lambda_x\in k$, is equal to 0. For all $x\in\nolinebreak{[G_0]}$, choose $\widetilde{\lambda}_x\in \CO$ such that $\rho(\widetilde{\lambda}_x)=\lambda_x$. By~(\ref{v viz gamma}), we get an element $\sum_{x\in[G_0]}\limits \widetilde{\lambda}_x\frac{v_x}{|C_G(x)|}$ of $\CO\Proj(kG)$ whose modular character has values in the maximal ideal $J(\CO)$ of $\CO$. But by~(\ref{character}), the value at $g\in G_{p'}$ of this modular character is equal to
$$\sum_{x\in[G_0]}\limits \widetilde{\lambda}_x\frac{\Phi_{v_x}(g)}{|C_G(x)|},$$
which is equal to 0 if $g\notin G_0$, and to $\widetilde{\lambda}_x$ if $g$ is conjugate to $x\in [G_0]$ in $G$. It follows that $\widetilde{\lambda}_x\in J(\CO)$, hence $\lambda_x=\rho(\widetilde{\lambda}_x)=0$. Since $g\in G_{p'}$ was arbitrary, we get that $\lambda_x=0$ for any $x\in [G_0]$, so the elements $\gamma_x$, for $x\in [G_0]$, are linearly independent. This completes the proof of Proposition~\ref{basis}.\endpf
\section{The simple functors $S_P$}
In this section, we generalize the results of Section~\ref{S un} to the functor $S_P$, for an arbitrary finite $p$-group $P$.
\begin{enonce}{Theorem} \label{generators kP}Let $P$ be a finite $p$-group. Then for any finite group $G$, the evaluation $\genk{P}(G)$ is generated by the elements of $kR_k(G)$ of the form
$$\Ind_{RC_G(R)}^G\Inf_{RC_G(R)/R}^{RC_G(R)}k\sfc^{RC_G(R)/R}(F),$$
where $R$ is a subgroup of $G$ such that $R\hkr P$, and $F$ is an indecomposable projective $kRC_G(R)/R$-module.
\end{enonce}
\pf By definition $\genk{P}(G)=kT^\Delta(G,P)(k_P)$, and $T^\Delta(G,P)$ is generated by the bimodules of the form $\Ind_{N_{R,\pi,Q}}^{G\times P}\Inf_{\sur{N}_{R,\pi,Q}}^{N_{R,\pi,Q}}E$, where:
\begin{itemize}
\item $N_{R,\pi,Q}$ is the normalizer in $G\times P$ of a diagonal $p$-subgroup $\Delta(R,\pi,Q)=\{(\pi(q),q)\mid q\in Q\}$, where $\pi:Q\to R$ is a group isomorphism from a subgroup $Q$ of $P$ to a subgroup $R$ of $G$.
\item $\sur{N}_{R,\pi,Q}=N_{R,\pi,Q}/\Delta(R,\pi,Q)$. 
\item $E$ is an indecomposable projective $k\sur{N}_{R,\pi,Q}$-module.
\end{itemize}
Then $N_{R,\pi,Q}=\big\{(a,b)\in G\times P\mid \forall q\in Q,\;{^a\pi(q)}=\pi(^bq)\big\}$. Let $\hat{Q}\leq P$ denote the second projection of $N_{R,\pi,Q}$. Then $\sur{N}_{R,\pi,Q}$ fits in a short exact sequence
$$\xymatrix{
1\ar[r]&C_G(R)\ar[r]^-i&\sur{N}_{R,\pi,Q}\ar[r]^-s&\hat{Q}/Q\ar[r]&1,
}
$$
where $i$ is the map $x\mapsto (x,1)\Delta(R,\pi,Q)$ from $C_G(R)$ to $\sur{N}_{R,\pi,Q}$, and $s$ maps $(a,b)\Delta(R,\pi,Q)\in\sur{N}_{R,\pi,Q}$ to $bQ\in\hat{Q}/Q$. Since $\hat{Q}/Q$ is a $p$-group and $k$ is algebraically closed, it follows from Corollary 5.12.4 of~\cite{linckelmann-I} that there exists an indecomposable projective $kC_G(R)$-module $M$ such that
$$E\cong \Ind_{C_G(R)}^{\sur{N}_{R,\pi,Q}}M.$$
Then
\begin{align*}
\Ind_{N_{R,\pi,Q}}^{G\times P}\Inf_{\sur{N}_{R,\pi,Q}}^{N_{R,\pi,Q}}E&\cong \Ind_{N_{R,\pi,Q}}^{G\times P}\Inf_{\sur{N}_{R,\pi,Q}}^{N_{R,\pi,Q}}\Ind_{C_G(R)}^{\sur{N}_{R,\pi,Q}}M\\
&\cong \Ind_{N_{R,\pi,Q}}^{G\times P}\Ind_{\Delta(R,\pi,Q)(C_G(R)\times 1)}^{N_{R,\pi,Q}}\Inf_{C_G(R)}^{\Delta(R,\pi,Q)(C_G(R)\times 1)}M\\
&\cong \Ind_{\Delta(R,\pi,Q)(C_G(R)\times 1)}^{G\times P}\Inf_{C_G(R)}^{\Delta(R,\pi,Q)(C_G(R)\times 1)}M.
\end{align*}
Now  the tensor product $T:=\left(\Ind_{N_{R,\pi,Q}}^{G\times P}\Inf_{\sur{N}_{R,\pi,Q}}^{N_{R,\pi,Q}}E\right)\otimes_{kP}k_P$ can be viewed as
$$\left(\Ind_{X}^{G\times P}\Inf_{C_G(R)}^{X}M\right)\otimes_{kP}\left(\Ind_{Y}^{P\times 1}k_P\right),$$
where $X=\Delta(R,\pi,Q)(C_G(R)\times 1)\leq G\times P$ and $Y=P\times 1\leq P\times 1$.  It follows from~\cite{tensorbimodules} that this tensor product is isomorphic to 
\begin{moneq}\label{tensor}
T\cong\mathop{\bigoplus}_{x\in p_2(X)\dom P/p_1(Y)}\limits \Ind_{X*{^{(x,1)}Y}}^{G\times 1}\Big(\big(\Inf_{C_G(R)}^{X}M\big)\otimes_{k_2(X)\cap {^xk}_1(Y)}k_P\Big),
\end{moneq}
{\flushleft where} $p_2(X)$ is the second projection of $X$, and $p_1(Y)$ the first projection of $Y=P\times 1$. In particular $p_1(Y)=P$, so $p_2(X)\dom P/p_1(Y)$ consists of a single double coset, and we can assume $x=1$ in~(\ref{tensor}). \par
Moreover $k_1(Y)=\{y\in P\mid (y,1)\in Y\}=P$, and 
\begin{align*}
k_2(X)&=\{y\in P\mid (1,y)\in X\}=\{q\in Q\mid \pi(q)\in C_G(R)\}\\
&=\pi^{-1}\big(R\cap C_G(R)\big)=Z(Q)\cong Z(R).
\end{align*}
Now $X*Y=\{(a,b)\in G\times 1\mid \exists c\in P,\;(a,c)\in X,\,(x,1)\in Y\}$ is equal to the first projection of $X$, that is $RC_G(R)$. We get finally
\begin{align*}
T&\cong \Ind_{RC_G(R)\times 1}^{G\times 1}\left(\big(\Inf_{C_G(R)}^{\Delta(R,\pi,Q)(C_G(R)\times 1)}M\big)\otimes_{kZ(Q)}k_P\right)\\
&\cong\Ind_{RC_G(R)}^G\Inf_{C_G(R)/Z(R)}^{RC_G(R)}M_{Z(R)},
\end{align*}
where $M_{Z(R)}=M\otimes_{kZ(R)}k$, viewed as a $kC_G(R)/Z(R)$-module. The inflation symbol $\Inf_{C_G(R)/Z(R)}^{RC_G(R)}$ stands more precisely for $\Inf_{RG_G(R)/R}^{RC_G(R)}\Iso_{C_G(R)/Z(R)}^{RC_G(R)/R}$. \par 
To complete the proof of Theorem~\ref{generators kP}, it remains to observe that the construction $M\mapsto M_{Z(R)}$ induces a bijection between the set of isomorphism classes of projective indecomposable $kC_G(R)$-modules and the set of isomorphism classes of projective indecomposable $kC_G(R)/Z(R)$-modules. Transporting this module via the isomorphism $C_G(R)/Z(R)\to RC_G(R)/R$ gives the indecomposable projective $kRC_G(R)/R$-module $F:=\Iso_{C_G(R)/Z(R)}^{RC_G(R)/R}M_{Z(R)}$. Moreover $F$ has to be viewed as an element of $kR_k\big(RC_G(R)/R\big)$, that is as $k\sfc^{RC_G(R)/R}(F)$. This completes the proof of Theorem~\ref{generators kP}.\endpf
\begin{enonce}{Definition} Let $G$ be a finite group, and $x$ be a $p$-regular element of $G$.
\begin{enumerate}
\item Let $P$ be a $p$-subgroup of $G$. Then $x$ is said to {\em have defect $P$} in $G$ if $P$ is conjugate in $G$ to a Sylow $p$-subgroup of $C_G(x)$. A conjugacy class of $p$-regular elements of $G$ is said to have defect $P$ if some of its elements (or equivalently all of its elements) have defect $P$.
\item Let $P$ be a finite $p$-group. Then $x$ is said to have {\em defect isomorphic to $P$} if $P$ is isomorphic to a Sylow $p$-subgroup of $C_G(x)$. A conjugacy class of $p$-regular elements of $G$ is said to have defect isomorphic to $P$ if some of its elements (or equivalently all of its elements) have defect isomorphic to $P$.
\end{enumerate}
\end{enonce}
\begin{rem}{Remark}
Note that with this definition, the elements of {\em defect zero} of~$G$ are the elements with {\em trivial defect} (that is, equal to $\{1\}$) in $G$. This slight discordance in terminology is the same as in block theory, where the blocks of defect 0 are also the blocks with trivial defect.
\end{rem}
\begin{enonce}{Lemma} \label{defect R}Let $G$ be a finite group, and $R$ be a $p$-subgroup of $G$.
\begin{enumerate}
\item Let $xR\in \big(RC_G(R)/R\big)_{p'}$. Then $xR$ has defect 0 in $RC_G(R)/R$ if and only if there exists an element $y\in C_G(R)_{p'}\cap xR$ with defect $R$ in  $RC_G(R)$, i.e. such that $RC_G(R,y)/R$ is a $p'$-group. 
\item Let $y\in C_G(R)_{p'}$ such that $yR\in \big(RC_G(R)/R\big)_0$. Then 
$$\Ind_{RC_G(R)}^G\Inf_{RC_G(R)/R}^{RC_G(R)}\gamma_{RC_G(R)/R,\,yR}=0\in kR_k(G)$$
unless $y$ has defect $R$ in $G$.
\end{enumerate}
\end{enonce}
\pf 1. It is well known that since $R$ is a $p$-group, a $p'$ element $xR$ of $RC_G(R)/R$ can be lifted to a $p$'-element $y\in xR$ of $RC_G(R)$.  Moreover $\big(RC_G(R)\big)_{p'}=C_G(R)_{p'}$, so $y\in C_G(R)_{p'}$, and the centralizer $H$ of $yR=xR$ in $RC_G(R)/R$ is the image in $RC_G(R)/R$ of the centralizer of $y$ in $RC_G(R)$. Thus $H=RC_G(R,y)/R$, and $H$ is a $p'$-group if and only if $R$ is a Sylow subgroup of $RC_G(R,y)$, that is if and only if $y$ has defect $R$ in $RC_G(R)$.\spn
2. The normalizer $N_G(R)$ normalizes $RC_G(R)$, and it acts on the conjugacy classes of $RC_G(R)/R$ by permutation. Let $T_y$ be the stabilizer in $N_G(R)$ of the conjugacy class of $yR$ in $RC_G(R)/R$. Then $g\in T_y$ if and only if $^gy\in {^c}(yR)={^cy}R$ for some $c\in RC_G(R)$, that is $^gy={^cy}r$, for some element~$r$ of $R$. But $y\in C_G(R)$, so $^cy\in C_G(^cR)=C_G(R)$, hence the $p'$-element $^cy$ commutes with the $p$-element $r$. Then $r$ is the $p$-part of ${^cy}r={^gy}$ which is a $p'$-element. This forces $r=1$, and $^gy={^cy}$. Then $g\in cC_G(y)$, hence $g\in RC_G(R)C_G(y)\cap N_G(R)=RC_G(R)N_G(R,y)$. Thus $T_y\leq RC_G(R)N_G(R,y)$. Conversely, it is clear that $RC_G(R)N_G(R,y)$ stabilizes the conjugacy class of $yR$ in $RC_G(R)/R$, so $T_y=RC_G(R)N_G(R,y)$.\par
Let $S$ be a Sylow $p$-subgroup of $T_y$ containing $R$. For simplicity, we set $W_y:=\Inf_{RC_G(R)/R}^{RC_G(R)}\Gamma_{RC_G(R)/R,\,yR}$ and $w_y:=\Inf_{RC_G(R)/R}^{RC_G(R)}\gamma_{RC_G(R)/R,\,yR}$. Then
$$\Ind_{RC_G(R)}^Gw_y=\Ind_{SC_G(R)}^G\Ind_{RC_G(R)}^{SC_G(R)}w_y,$$
and we want to show that this is equal to zero unless $y$ has defect $R$ in $G$. We will show a little more: We claim that $\Ind_{RC_G(R)}^{SC_G(R)}w_y=0$ in $kR_k\big(SC_G(R)\big)$ if $R$ is not a Sylow $p$-subgroup of $C_G(y)$. \par
To prove this claim, we compute the coefficients of $\Ind_{RC_G(R)}^{SC_G(R)}w_y$ in the basis of $kR_k\big(SC_G(R)\big)$ consisting of the simple $k\big(SC_G(R)\big)$-modules. Since $w_y$ is the reduction in $k$ of $W_y$, we can use the modular character $\theta_y$ associated to $W_y$ to do this computation. Let $U$ be a simple $k\big(SC_G(R)\big)$-module. Since $RC_G(R)$ is a normal subgroup of $SC_G(R)$ with $p$-power index, it follows from Corollary 5.12.4 of~\cite{linckelmann-I} that the projective cover of $U$ is isomorphic to $\Ind_{RC_G(R)}^{SC_G(R)}E$, where $E$ is a projective $kRC_G(R)$-module. The coefficient of $U$ in the expression of $\Ind_{RC_G(R)}^{SC_G(R)}w_y$ in the basis of simple $kSC_G(R)$-modules is then equal to $\rho(m_U)$, where
\begin{align*}
m_U=\langle \Ind_{RC_G(R)}^{SC_G(R)}E,\Ind_{RC_G(R)}^{SC_G(R)}\theta_y\rangle_{SC_G(R)}&=\langle E,\Res_{RC_G(R)}^{SC_G(R)}\Ind_{RC_G(R)}^{SC_G(R)}\theta_y\rangle_{RC_G(R)}\\
&=\sum_{g\in SC_G(R)/RC_G(R)}\langle E,{^g\theta_y}\rangle_{RC_G(R)}\\
&=|SC_G(R)/RC_G(R)|\langle E,\theta_y\rangle_{RC_G(R)},
\end{align*}
since $^g\theta_y=\theta_y$ for $g\in SC_G(R)$, as $SC_G(R)\leq T_y$. \par
Then $\rho(m_U)=0$ in $k$ if $SC_G(R)\neq RC_G(R)$, or equivalently if $p$ divides $|T_y:RC_G(R)|=|N_G(R,y):RC_G(R,y)|$, that is, since $RC_G(R,y)/R$ is a $p'$-group by Assertion~1, if $p$ divides $|N_G(R,y):R|$.  Hence if $\rho(m_U)\neq 0$, then $R$ is a Sylow $p$-subgroup of $N_G(R,y)=N_{C_G(y)}(R)$, hence $R$ is a Sylow $p$-subgroup of $C_G(y)$, which proves the claim. This completes the proof of Lemma~\ref{defect R}.\endpf
\begin{enonce}{Corollary} \label{defect kP} Let $G$ be a finite group, and $P$ be a finite $p$-group. Let $[G_{\hkr P}]$ be a set of representatives of conjugacy classes of the set $G_{\hkr P}$ of $p$-regular elements of $G$ with defect $R\hkr P$. For $x\in G_{\hkr P}$, let $R_x$ be a chosen Sylow $p$-subgroup of $C_G(x)$. Then the elements
$$U_x=\Ind_{R_xC_G(R_x)}^G\Inf_{R_xC_G(R_x)/R_x}^{R_xC_G(R_x)}\gamma_{R_xC_G(R_x)/R_x,\, xR_x}$$
for $x\in [G_{\hkr P}]$, form a basis of $\genk{P}(G)$.
\end{enonce}
\pf By Lemma~\ref{defect R}, the elements $U_x$, for $x\in [G_{\hkr P}]$, generate $\genk{P}(G)$, and all we have to show is that these elements are linearly independent. Let $S$ be a Sylow $p$-subgroup of $G$. If $x\in [G_{\hkr P}]$ has defect $R\leq G$, then $R\hkr S$, so $U_x\in\genk{S}(G)$. So it is enough to prove that the elements $U_x$, for $x\in [G_{\hkr S}]$, are linearly independent. But $\genk{S}(G)=kR_k(G)$ by Corollary~\ref{Sylow}, and $\big|[G_{\hkr S}]\big|=\big|[G_{p'}]\big|$ as any element of $G_{p'}$ as defect $R$ for some $R\hkr S$. Hence
$$\big|[G_{p'}]\big|=\dim_k kR_k(G)=\dim_k\genk{S}(G)\leq \big|[G_{\hkr S}]\big|=\big|[G_{p'}]\big|,$$
so $\dim_k\genk{S}(G)= \big|[G_{\hkr S}]\big|$, which completes the proof. \endpf
\begin{rem}{Remark} Observe that for $x\in G_{\hkr P}$, the element $U_x$ does not depend on the choice of the Sylow $p$-subgroup $R_x$. Moreover $U_{x}=U_{x^g}$, for any $g\in G$. So the set $\{U_x\mid x\in [G_{\hkr P}]\}$ is a canonical basis of $\genk{P}(G)$.
\end{rem}
\begin{enonce}{Theorem} Let $G$ be a finite group, and $P$ be a finite $p$-group. Let $G_{\cong P}$ denote the set of elements of $G$ with defect isomorphic to $P$, and $[G_{\cong P}]$ be a set of representatives of conjugacy classes of elements of $G_{\cong P}$. Then the images of the elements $U_x$, for $x\in [G_{\cong P}]$, under the projection $\genk{P}(G)\twoheadrightarrow S_P(G)$, form a basis of $S_P(G)$.
\end{enonce}
\pf Let $x\in G_{\hkr P}$, and $R=R_x$ be a Sylow subgroup of $C_G(x)$. Then the element $U_x$ defined in Corollary~\ref{defect kP} lies in $\genk{R}(G)$, so $U_x$ is sent to 0 under the projection $\genk{P}(G)\twoheadrightarrow S_P(G)$ if $R\ncong P$. It follows that the images of the elements $U_x$, for $x\in [G_{\cong P}]$, under the projection $\genk{P}(G)\twoheadrightarrow S_P(G)$, generate $S_P(G)$. In particular $\dim_kS_P(G)\leq \big|[G_{\cong P}]\big|$.\par
Now, as in the proof of Theorem~\ref{composition}, we choose an enumeration $P_1, P_2,\ldots$ of $\cal P$ with the property that $P_i\hkr P_j$ implies $i\leq j$. This gives a filtration
$$0=F_0< F_1<\ldots <F_n< F_{n+1}<\ldots$$
of $kR_k$ by subfunctors $F_i$, for $i\in\NN$, such that $F_i/F_{i-1}\cong S_{P_i}$ for $i>0$. \par
Let $n$ be the unique integer such that $P_n$ is isomorphic to a Sylow $p$-subgroup of $G$. If $i>n$, then $S_{P_i}(G)=0$ by Lemma~\ref{SP},  since $P_i$ is not isomorphic to a subgroup of $G$: Indeed, if it were, then $P_i$ would be isomorphic to a subgroup of $P_n$, which would imply $i\leq n$. It follows that
$$F_n(G)=F_{n+1}(G)=\ldots=kR_k(G).$$
So we have a filtration of $kR_k(G)$
$$0\leq F_1(G)\leq\ldots\leq F_{i-1}(G)\leq F_i(G)\leq\ldots\leq F_n(G)=kR_k(G).$$
Moreover, if the quotient $F_i(G)/F_{i-1}(G)\cong S_{P_i}(G)$ is non zero, then $P_i\hkr G$ by Lemma~\ref{SP}, i.e. $P_i\hkr P_n$, so $i\leq n$. Then
\begin{align*}
\big|[G_{p'}]\big|&=\dim_k kR_k(G)=\sum_{i=1}^n\dim_kS_{P_i}(G)\\
&=\sum_{R\hkr P}\dim_k S_R(G)\leq \sum_{R\hkr P}\big|[G_{\cong R}]\big|=\big|[G_{p'}]\big|.
\end{align*}
Hence all inequalities $\dim_k S_R(G)\leq \big|[G_{\cong R}]\big|$ are equalities, and  the theorem follows.\spn
\begin{enonce}{Corollary} Let $G$ be a finite group, and $P$ be a finite $p$-group. Then the dimension of $S_P(G)$ is equal to the number of conjugacy classes of $p$-regular elements of $G$ with defect isomorphic to $P$.
\end{enonce}
\begin{rem}{Remark} This corollary is consistent with Corollary~6.14 of~\cite{diagonal-functors-in-char-p}, via Remark~\ref{other parametrization}.
\end{rem}
\vspace{2ex}\mpn
{\bf Acknowledgements :} I wish to thank Deniz Y\i lmaz for his careful reading of a preliminary version of this work, and his pertinent remarks and suggestions.

\centerline{\rule{5ex}{.1ex}}
\begin{flushleft}
Serge Bouc, CNRS-LAMFA, Universit\'e de Picardie, 33 rue St Leu, 80039, Amiens, France.\\
{\tt serge.bouc@u-picardie.fr}\\
\end{flushleft}
\end{document}